%

\input xy
\xyoption{all}



\font\sc=cmcsc10 \rm


\newcount\secnb
\newcount\parnb
\secnb=0
\parnb=0
\newtoks\secref
\secref={}

\def\smallskip{\par\vskip 2mm}
\def\medskip{\par\vskip 5mm}
\def\goodbreak{\penalty -100}

\def\section#1{\global\advance\secnb by 1
\secref=\expandafter{\the\secnb.}\parnb=0
\medskip\goodbreak
\centerline{\bf\S\the\secref\ #1}
\smallskip\nobreak}

\def\references{\medskip\goodbreak
\centerline{\bf R\'ef\'erences}
\smallskip\nobreak}

\def\tit{\global\advance\parnb by 1
\smallskip\goodbreak
\noindent\the\secref\the\parnb)\ }

\def\endtit{\par\nobreak}



\def\label#1{\relax}

\def\ref#1{\csname crossref#1\endcsname}


\newcount\bibrefnb
\bibrefnb=0

\def\biblabel#1{\global\advance\bibrefnb by 1\item{\the\bibrefnb.}}

\def\refto#1{\setbox1=\hbox{\sc #1}\copy1}
\def\reftosame{\hbox to \the\wd1{\hrulefill}}

\def\bibref#1{\csname bibref#1\endcsname}

\def\bye{\medskip\leftline{{\it Addresse:}}{\leftskip=1cm\parindent=0cm

Laboratoire J.A. Dieudonn\'e,\hfill\break
Universit\'e de Nice,\hfill\break
Parc Valrose,\hfill\break
F-06108 Nice Cedex 02 (France).\par}

\smallskip\leftline{{\it Courriel:}}{\leftskip=1cm\parindent=0cm

Clemens Berger $<$cberger@math.unice.fr$>$,\hfill\break
Benoit Fresse $<$fresse@math.unice.fr$>$.\par}

\vfill\supereject\end}

\def\beginresume{\smallskip
\begingroup\leftskip=1cm\rightskip=1cm
{\sc R\'esum\'e}. --- }
\def\endresume{\endgroup}
\def\beginabstract{\smallskip
\begingroup\leftskip=1cm\rightskip=1cm
{\sc Abstract}. --- }
\def\endabstract{\endgroup}

\def\W{{\cal W}}
\def\E{{\cal E}}
\def\M{{\cal M}}
\def\F{{\cal F}}
\def\X{{\cal X}}
\def\Z{{\bf Z}}
\def\N{{\bf N}}
\def\TR{\mathop{TR}\nolimits}
\def\TC{\mathop{TC}\nolimits}
\def\Id{\mathop{Id}\nolimits}


\hsize=15cm
\vsize=20.5cm

\headline={\ifnum\pageno=1{\it Projet de note aux C.R.Acad.Sci.Paris S\'er. I Math. (Topologie)}\hfil\else\hfil\fi}

\vbox{\medskip
\centerline{\bf Une d\'ecomposition prismatique de l'op\'erade de Barratt-Eccles}
\smallskip
\centerline{\sc Clemens Berger, Benoit Fresse}
\smallskip
\centerline{26/4/2002}
\medskip
\beginresume
L'op\'erade de Barratt-Eccles est une op\'erade simpliciale form\'ee par les constructions bar homog\`enes des groupes sym\'etriques.
On montre que ces ensembles simpliciaux se d\'ecomposent en r\'eunions de prismes index\'es par des surjections.
On observe que les complexes cellulaires d\'efinis par cette structure prismatique
s'identifient aux composantes de l'op\'erade des surjections
(l'op\'erade introduite par J. McClure et J. Smith dans leur travaux sur la conjecture de Deligne).
\endresume
\medskip
\centerline{\bf A prismatic decomposition of the Barratt-Eccles operad}
\smallskip
\beginabstract
The Barratt-Eccles operad is a simplicial operad formed by the classical homogeneous bar construction of the symmetric groups.
We prove that these simplicial sets decompose as unions of prisms indexed by surjections.
We observe that the cellular complexes given by this prismatic structure
are nothing but the components of the surjection operad
(the operad introduced by J. McClure and J. Smith in their work on the Deligne conjecture).
\endabstract
\smallskip}

\medskip
\noindent{\bf Abridged English Version}
\medskip

The purpose of this note is to discuss the relationship between the {\it Barratt-Eccles operad} and the {\it surjection operad}.
We prove that the surjection operad arises from a prismatic decomposition of the Barratt-Eccles operad.
This work sheds a new light on results of J. McClure and J. Smith
(cf. [\bibref{MCSPrism}], [\bibref{MCSSeq}]).

\tit{\it The Barratt-Eccles operad} ---
We follow the conventions of our former article (cf. [\bibref{BF}]).
The simplicial Barratt-Eccles operad is denoted by $\W$.
The letter $\E$ denotes the normalized differential graded operad associated to $\W$.
The simplicial set $\W(r)$ is the homogeneous bar construction of the symmetric group $\Sigma_r$.
Hence, an $n$-dimensional simplex in $\W(r)$ is an $n+1$-tuple of permutations $(w_0,\ldots,w_n)\in\W(r)_n$.
We have $d_i(w_0,\ldots,w_n) = (w_0,\ldots,\widehat{w_i},\ldots,w_n)$
and $s_j(w_0,\ldots,w_n) = (w_0,\ldots,w_j,w_j,\ldots,w_n)$.
The differential graded module $\E(r)$ is the normalized chain complex of $\W(r)$.
We refer to the literature for more details about the operad structure of $\W(r)$ and $\E(r)$.

The surjection operad is denoted by the letter $\X$.
We recall that the module $\X(r)_d$ is generated by the {\it non-degenerated surjections}
$u: \{1,\ldots,r+d\}\,\rightarrow\,\{1,\ldots,r\}$.
A surjection is non-degenerated if $u(i+1)\not=u(i)$, for $i = 1,\ldots,r+d-1$.

Let us mention that the degree $0$ component of $\E(r)$ and $\X(r)$
is nothing but the regular representation of the symmetric group:
$\E(r)_0 = \X(r)_0 = \Z[\Sigma_r]$.

\tit{\it Cellular $E$-infinity structures} ---
The Barratt-Eccles operad has a filtration 
$F_1\W\subset F_2\W\subset\cdots\subset F_n\W\subset\cdots\subset F_\infty\W = \W$
by simplicial operads $F_n\W$ whose topological realization is equivalent to the {\it little $n$-cubes operad} (cf. [\bibref{B}]).
We have a cellular structure which refines this filtration
and which we denote as the cellular $E$-infinity structure of the Barratt-Eccles operad.
The cells $F_{n_{i j},w}\W(r)\subset\W(r)$ are indexed
by pairs
consisting of a collection of non-negative integers $n_{i j}\in\N$ and of a permutation $w\in\Sigma_r$.
We have explicitly $F_n\W(r) = \bigcup_{n_{i j}<n} F_{n_{i j},w}\W(r)$.
We have an induced cellular structure and an induced filtration on the differential graded operad $\E$.

The surjection operad is also equipped with a cellular $E$-infinity structure.
The reader is refered to the article of J. McClure and J. Smith (cf. [\bibref{MCSSeq}]).

\tit{\sc Theorem} --- {\it There are chain-morphisms $\TR: \E\,\rightarrow\,\X$ and $\TC: \X\,\rightarrow\,\E$
such that:

(1) in degree $0$, the morphisms $\TR$ and $\TC$ are the identity of $\Z[\Sigma_r]$;

(2) we have $\TR\cdot\TC = \Id_{\X}$;

(3) the morphisms $\TR$ and $\TC$ preserve the cellular $E$-infinity structures and the operad filtrations;

(4) the map $\TR$ is an operad morphism.}

\smallskip
The map $\TR: \E\,\rightarrow\,\X$ is the morphism introduced in the article [\bibref{BF}].
Assertion (4) is the main result of this article.
The purpose of this note is to construct the section $\TC: \X\,\rightarrow\,\E$ of $\TR: \E\,\rightarrow\,\X$.

\tit{\sc Lemma} --- {\it There is a map $H: \E(r)_*\,\rightarrow\,\E(r)_{*+1}$
such that $\TC\cdot\TR = \Id_{\E} + H\cdot\delta + \delta\cdot H$.
This map preserves the cellular $E$-infinity structure of $\E$.}

\smallskip
This lemma is a consequence of a general property.
Let $T: \E(r)_*\,\rightarrow\,\E(r)_*$ be a chain map which is the identity morphism in degree $0$.
We have in fact a canonical homotopy $H: \E(r)_*\,\rightarrow\,\E(r)_{*+1}$
defined by $H(w_0,\ldots,w_d) = \sum_{i=0}^d (-1)^i (T(w_0,\ldots,w_i),w_i,\ldots,w_d)$.
To be more precise, the expression $T(w_0,\ldots,w_i)$ represents a sum of $i$-dimensional simplices in $\W(r)$.
These are $i+1$-tuples of permutations which can be concatenated with $(w_i,\ldots,w_d)$.
This process gives a sum of $d+1$-dimensional simplices in $\W(r)$
and hence an element of $\E(r)_{*+1}$.
The relation $T = \Id_{\E} + H\cdot\delta + \delta\cdot H$ is readily verified.
In addition, we observe that $H$ preserves the cellular structure of $\E$
if $T$ satisfies this property.

\smallskip
The next theorem is an immediate corollary of the results above:

\tit{\sc Theorem} --- {\it The morphism $\TR: \E\,\rightarrow\,\X$ induces a quasi-isomorphism of operads
$F_n\TR: F_n\E\,\rightarrow\,F_n\X$, for $n = 1,2,\ldots,\infty$.}

\medskip
\hbox to 3cm{\hrulefill}

\medskip
Le but de cette note est de pr\'eciser la relation entre deux op\'erades $E$-infinis:
l'{\it op\'erade de Barratt-Eccles}, d'une part,
l'{\it op\'erade des surjections}, d'autre part.
On compl\`ete les r\'esultats obtenus dans l'article [\bibref{BF}].
On montre que l'op\'erade des surjections provient d'une d\'ecomposition prismatique de l'op\'erade de Barratt-Eccles.
Ce travail donne un nouveau point de vue sur des r\'esultats de J. McClure et J. Smith
(cf. [\bibref{MCSPrism}], [\bibref{MCSSeq}]).

\section{R\'esultats}

\tit{\it L'op\'erade de Barratt-Eccles} ---
On reprend les conventions de notre article (cf. [\bibref{BF}]).
L'op\'erade de Barratt-Eccles simpliciale est not\'ee $\W$.
L'op\'erade diff\'erentielle gradu\'ee associ\'ee est d\'esign\'ee par la lettre $\E$.
L'ensemble simplicial $\W(r)$ est la construction bar homog\`ene du groupe sym\'etrique $\Sigma_r$.
Ainsi, les simplexes de dimension $n$ de $\W(r)$
sont les $n+1$-uplets de permutations $(w_0,\ldots,w_n)\in\W(r)_n$.
On a $d_i(w_0,\ldots,w_n) = (w_0,\ldots,\widehat{w_i},\ldots,w_n)$
et $s_j(w_0,\ldots,w_n) = (w_0,\ldots,w_j,w_j,\ldots,w_n)$.
Le module diff\'erentiel gradu\'e $\E(r)$ est le complexe des cha\^\i nes normalis\'ees de $\W(r)$.
On renvoie le lecteur \`a la litt\'erature pour plus de d\'etails sur les structures d'op\'erades de $\W(r)$ et de $\E(r)$.

L'op\'erade des surjections est d\'esign\'ee par la lettre $\X$.
On rappelle que le module $\X(r)_d$ est engendr\'e par les surjections non-d\'eg\'en\'er\'ees
$u: \{1,\ldots,r+d\}\,\rightarrow\,\{1,\ldots,r\}$.
On dit qu'une surjection $u$ est non-d\'eg\'en\'er\'ee si on a $u(i+1)\not=u(i)$, pour $i = 1,\ldots,r+d-1$.

On note que la composante de degr\'e $0$ de $\E(r)$ et de $\X(r)$
est la repr\'esentation r\'eguli\`ere du groupe sym\'etrique
$\E(r)_0 = \X(r)_0 = \Z[\Sigma_r]$.

\tit{\it Structures cellulaires $E$-infini} ---
On a une filtration de l'op\'erade de Barratt-Eccles
$F_1\W\subset F_2\W\subset\cdots\subset F_n\W\subset\cdots\subset F_\infty\W = \W$
par des op\'erades simpliciales $F_n\W$
dont la r\'ealisation topologique est \'equivalente
\`a l'op\'erade des {\it petits $n$-cubes} (cf. [\bibref{B}]).
Cette filtration se raffine en une structure cellulaire
que l'on d\'esignera comme la {\it structure cellulaire $E$-infini}
de l'op\'erade de Barratt-Eccles.
Les cellules $F_{n_{i j},w}\W(r)\subset\W(r)$ sont index\'ees par les paires constitu\'ees
d'une collection d'entiers $n_{i j}\in\N$
et d'une permutation $w\in\Sigma_r$.
On a explicitement $F_n\W(r) = \bigcup_{n_{i j}<n} F_{n_{i j},w}\W(r)$.
On consid\`ere la structure cellulaire $E$-infini induite sur l'op\'erade diff\'erentielle gradu\'ee $\E$.

\smallskip
Le but de cette note est d'\'etablir le th\'eor\`eme suivant:

\tit{\sc Th\'eor\`eme} --- {\it On a des morphismes de complexes de cha\^\i nes
$\TR: \E\,\rightarrow\,\X$ et $\TC: \X\,\rightarrow\,\E$
tels que:

(1) en degr\'e $0$, les morphismes $\TR$ et $\TC$ sont l'identit\'e de $\Z[\Sigma_r]$;

(2) on a $\TR\cdot\TC = \Id_{\X}$;

(3) le morphisme compos\'e $\TC\cdot\TR: \E\,\rightarrow\,\E$
pr\'eserve la structure cellulaire $E$-infini de $\E$;

(4) l'application $\TR$ est un morphisme d'op\'erades.}

\smallskip
L'application $\TR: \E\,\rightarrow\,\X$ est le morphisme introduit dans l'article [\bibref{BF}].
L'assertion (4) du th\'eor\`eme est le r\'esultat principal de cet article.
L'application $\TC: \X\,\rightarrow\,\E$ est d\'efinie dans la section 3 de cette note.
L'assertion (2) r\'esulte d'un \'enonc\'e plus pr\'ecis (cf. lemme \ref{AWPpty}).

On a une structure cellulaire $E$-infini sur l'op\'erade des surjections qui a \'et\'e introduite par J. McClure et J. Smith.
On prouve dans l'article [\bibref{BF}] que le module $F_{n_{i j},w}\X(r)$
d\'efini par ces auteurs
est l'image de $F_{n_{i j},w}\E(r)$ par le morphisme $\TR: \E\,\rightarrow\,\X$.
Cette propri\'et\'e se d\'eduit aussi de la compatibilit\'e des structures cellulaires avec les structures d'op\'erades.
Le lecteur v\'erifiera facilement en reprenant la d\'efinition des structures cellulaires et les arguments de l'article [\bibref{BF}]
que l'image de $F_{n_{i j},w}\X(r)$ par le morphisme $\TC: \X\,\rightarrow\,\E$
est contenue dans $F_{n_{i j},w}\E(r)$.
Ceci prouve l'assertion (3) du th\'eor\`eme.

\smallskip
On a aussi la propri\'et\'e suivante:

\tit{\sc Lemme} --- {\it On a une application $H: \E(r)_*\,\rightarrow\,\E(r)_{*+1}$
telle que $\TC\cdot\TR = \Id_{\E} + H\cdot\delta + \delta\cdot H$.
Cette application pr\'eserve la structure cellulaire $E$-infini de $\E$.}

\smallskip
\noindent{\it Preuve:} Ce lemme est cons\'equence d'un r\'esultat plus g\'en\'eral:
on se donne un morphisme de complexes $T: \E(r)_*\,\rightarrow\,\E(r)_*$
qui est l'identit\'e en degr\'e $0$.
On a alors une homotopie naturelle $H: \E(r)_*\,\rightarrow\,\E(r)_{*+1}$
entre $T$ et l'identit\'e de $\E$.
Explicitement,
on pose $H(w_0,\ldots,w_d) = \sum_{i=0}^d (-1)^i (T(w_0,\ldots,w_i),w_i,\ldots,w_d)$.
L'expression $T(w_0,\ldots,w_i)$ repr\'esente une somme de simplexes de dimension $i$ dans $\W(r)$.
Ce sont des $i+1$-uplets de permutations que l'on concat\`ene avec $(w_i,\ldots,w_d)$.
On obtient ainsi des simplexes de dimension $d+1$ dans $\W(r)$
et donc un \'el\'ement de $\E(r)_{*+1}$.
La relation $T = \Id_{\E} + H\cdot\delta + \delta\cdot H$ se v\'erifie sans difficult\'es.
On constate \'egalement que $H$ pr\'eserve la structure cellulaire $E$-infini de $\E$
si $T$ a cette propri\'et\'e.

\smallskip
Les r\'esultats ci-dessus ont pour corollaire imm\'ediat:

\tit{\sc Th\'eor\`eme}\label{OpFilt} --- {\it Le morphisme $\TR: \E\,\rightarrow\,\X$
induit un quasi-isomorphisme d'op\'erades $F_n\TR: F_n\E\,\rightarrow\,F_n\X$,
quelque soit $n = 1,2,\ldots,\infty$.}

\section{La d\'ecomposition prismatique de l'op\'erade de Barratt-Eccles}

\tit{\it Le prisme associ\'e \`a une surjection}\label{PrismDefn} --- On fixe une surjection $u\in\X(r)_d$.
On note $d_k$ le nombre d'occurences de $k\in\{1,\ldots,r\}$ dans la suite $(u(1),\ldots,u(r+d))$.
(On a de la sorte $d+r = d_1+\cdots+d_r$.)

On a un prisme $\tau_u: \Delta^{d_1-1}\times\cdots\times\Delta^{d_r-1}\,\rightarrow\,\W(r)$ associ\'e \`a $u$.
L'image d'un simplexe $\sigma\in(\Delta^{d_1-1}\times\cdots\times\Delta^{d_r-1})_n$ dans $\W(r)_n$
est d\'etermin\'ee par l'image de ses sommets $(\sigma(0),\ldots,\sigma(n))$.
On a explicitement $\tau_u(\sigma) = (\tau_u(\sigma(0)),\ldots,\tau_u(\sigma(n))$.
Un sommet du prisme $\Delta^{d_1-1}\times\cdots\times\Delta^{d_r-1}$
est sp\'ecifi\'e par un $r$-uplet d'entiers $(x_1,\ldots,x_r)$
tels que $0\leq x_k\leq d_k-1$.
Le sommet correspondant $\tau_u(x_1,\ldots,x_r)\in\W(r)_0$
est la permutation de $(1,\ldots,r)$
qui est d\'efinie par la sous suite de $(u(1),\ldots,u(r+d))$
form\'ee par les occurrences num\'ero $x_1+1,\ldots,x_r+1$
des valeurs $1,\ldots,r$.

On consid\`ere par exemple la surjection $(u(1),\ldots,u(5)) = (1,2,3,1,2)$.
Le prisme associ\'e $\tau_u: \Delta^1\times\Delta^1\times\Delta^0\,\rightarrow\,\W(3)$
se repr\'esente par la figure suivante:
$$\vcenter{\xymatrix{ (0,1,0)\ar[r] & (1,1,0) \\
(0,0,0)\ar[u]\ar[ur]\ar[r] & (1,0,0)\ar[u] \\ }}\qquad\mapsto\qquad\vcenter{\xymatrix{ (1,3,2)\ar[r] & (3,1,2) \\
(1,2,3)\ar[u]\ar[ur]\ar[r] & (2,3,1)\ar[u] \\ }}$$

Les r\'esultats ci-dessous montrent que ces prismes d\'efinissent une d\'ecomposition cellulaire de $\W(r)$.
La preuve du lemme \ref{PrismDec} (cit\'e en remarque) est omise.

\tit{\sc Lemme}\label{SurjFaces} ---
{\it Les faces d'un prisme $\tau_u$ sont les prismes $\tau_v$ associ\'es aux sous-suites $(v(1),\ldots,v(r+e))$ de $(u(1),\ldots,u(r+d))$.}

\smallskip
\noindent{\it Preuve:} On consid\`ere par exemple la composition de $\tau_u$ avec le morphisme de face
$$1\times\cdots\times d^x\times\cdots\times 1:
\Delta^{d_1-1}\times\cdots\times\Delta^{d_k-2}\times\cdots\times\Delta^{d_r-1}
\,\rightarrow\,\Delta^{d_1-1}\times\cdots\times\Delta^{d_k-1}\times\cdots\times\Delta^{d_r-1}.$$
Ce morphisme compos\'e s'identifie \`a un prisme
$\tau_v: \Delta^{d_1-1}\times\cdots\times\Delta^{d_k-2}\times\cdots\times\Delta^{d_r-1} \,\rightarrow\,\W(r)$.
La surjection $v$ s'obtient
en omettant la $x+1$-i\`eme occurrence de $k$
dans la suite $(u(1),\ldots,u(r+d))$.
On g\'en\'eralise facilement cette construction \`a tous les sous-prismes
$\Delta^{e_1-1}\times\cdots\times\Delta^{e_r-1}\subset\Delta^{d_1-1}\times\cdots\times\Delta^{d_r-1}$.

\tit{\sc Lemme}\label{PrismDec}\endtit
{\it (1) On a $\tau_v(\Delta^{e_1-1}\times\cdots\times\Delta^{e_r-1})
\subset\tau_u(\Delta^{d_1-1}\times\cdots\times\Delta^{d_r-1})$
si et seulement si $(v(1),\ldots,v(r+e))$ est une sous-suite de $(u(1),\ldots,u(r+d))$.
Le prisme $\tau_v$ est alors une face de $\tau_u$.

(2) Les images des prismes $\tau_u$, $u\in\X(r)$, recouvrent l'ensemble simplicial $\W(r)$.
Les prismes $\tau_u$, $u\in\X(r)$, s'intersectent selon des r\'eunions de faces.}

\tit{\it Le simplexe fondamental associ\'e \`a une surjection} ---
On reprend le prisme associ\'e \`a une surjection donn\'ee $u: \{1,\ldots,r+d\}\,\rightarrow\,\{1,\ldots,r\}$.
On sp\'ecifie un {\it simplexe fondamental} parmi les simplexes maximaux de $\Delta^{d_1-1}\times\cdots\times\Delta^{d_r-1}$.
Un simplexe maximal $\sigma\in(\Delta^{d_1-1}\times\cdots\times\Delta^{d_r-1})_d$
est d\'etermin\'e par une suite d'entiers $k_i\in\{1,\ldots,r\}$, $i = 0,\ldots,d-1$.
On consid\`ere la suite de sommets $(x^{(i)}_1,\ldots,x^{(i)}_r)\in(\Delta^{d_1-1}\times\cdots\times\Delta^{d_r-1})_0$
telle que $x^{(i+1)}_k = x^{(i)}_k+1$ si $k=k_i$ et $x^{(i+1)}_k = x^{(i)}_k$ sinon.
On a alors un et un seul simplexe maximal dans $\Delta^{d_1-1}\times\cdots\times\Delta^{d_r-1}$
dont les sommets sont les $(x^{(i)}_1,\ldots,x^{(i)}_r)$
(cf. [\bibref{GZ}, Section II.5]).
Le simplexe fondamental de $\tau_u$
est le simplexe maximal associ\'e \`a la suite $(k_0,\ldots,k_{d-1})$
form\'ee par les {\it c\'esures} de la surjection $u$.
On obtient cette suite en retirant de $(u(1),\ldots,u(r+d))$
la derni\`ere occurrence de chaque valeur $k\in\{1,\ldots,r\}$ (cf. [\bibref{BF}]).

Le simplexe fondamental associ\'e \`a la surjection $u = (1,2,3,1,2)$ du paragraphe \ref{PrismDefn}
a pour sommets $((1,2,3),(2,3,1),(3,1,2))$.

\tit{\it Remarques:}
On peut consid\'erer le complexe cellulaire associ\'e \`a notre d\'ecomposition prismatique de $\W(r)$.
Ce complexe est isomorphe au module $\X(r)$.
Le lemme \ref{EZPpty} montre essentiellement
que la diff\'erentielle de l'op\'erade des surjections, d\'efinie dans l'article [\bibref{BF}],
correspond \`a la diff\'erentielle cellulaire.

Outre les travaux de McClure et Smith (cf. [\bibref{MCSPrism}], [\bibref{MCSSeq}]),
l'op\'erade $F_2\X$ appara\^\i t (sous la notation $\M$) dans le travail de Kontsevich et Soibelman
sur la conjecture de Deligne (cf. [\bibref{KS}]).
La r\'eunion des images des prismes index\'es par les surjections de $F_2\X$
est la collection simpliciale $\F$ de l'article [\bibref{MCSPrism}].

\section{Les morphismes}

\tit{\it Le morphisme d'Eilenberg-Zilber} ---
On d\'efinit l'\'el\'ement $\TC(u)\in\E(r)_d$ associ\'e \`a une surjection $u\in\X(r)_d$.
On prend la somme altern\'ee des simplexes maximaux du prisme
$\tau_u: \Delta^{d_1-1}\times\cdots\times\Delta^{d_r-1}\,\rightarrow\,\W(r)_d$.
Un simplexe a un signe positif si son orientation naturelle concorde avec l'orientation du simplexe fondamental,
n\'egatif sinon.
Ce morphisme diff\`ere de l'application d'Eilenberg-Zilber classique par un signe.
Plus explicitement,
l'application d'Eilenberg-Zilber associe
au g\'en\'erateur de $N_{d_1-1}(\Delta^{d_1-1})\otimes\cdots\otimes N_{d_r-1}(\Delta^{d_r-1})$
une somme de simplexes dans $N_d(\Delta^{d_1-1}\times\cdots\times\Delta^{d_r-1})$.
On prend l'image de cette somme par $N_*(\tau_u)$.
C'est $\pm\TC(u)\in\E(r)_d$.
Le signe est d\'etermin\'e par l'orientation du simplexe fondamental.
Ainsi,
pour l'exemple du paragraphe \ref{PrismDefn},
on obtient $\TC(1,2,3,1,2) = ((1,2,3),(2,3,1),(3,1,2)) - ((1,2,3),(1,3,2),(3,1,2))$.

\tit{\sc Lemme}\label{EZPpty} --- {\it L'application $\TC: \X\,\rightarrow\,\E$ est un morphisme de complexes.}

\smallskip
\noindent{\it Preuve:}\ On identifie l'ensemble simplicial $\Delta^{d_k-1}$
au simplexe de dimension $d_k-1$ qui l'engendre.
Ainsi, le produit tensoriel $\Delta^{d_1-1}\otimes\cdots\otimes\Delta^{d_r-1}$
repr\'esente le g\'en\'erateur de $N_{d_1-1}(\Delta^{d_1-1})\otimes\cdots\otimes N_{d_r-1}(\Delta^{d_r-1})$.

L'application d'Eilenberg-Zilber classique est un morphisme de complexes.
La somme des faces de $\Delta^{d_1-1}\otimes\cdots\otimes\Delta^{d_r-1}$
dans $N_*(\Delta^{d_1-1})\otimes\cdots\otimes N_*(\Delta^{d_r-1})$
est donc envoy\'ee sur la diff\'erentielle de $\TC(u)$ dans $N_*(\W(r))$.
L'image de la face
$\Delta^{d_1-1}\otimes\cdots\otimes d_x(\Delta^{d_k-1})\otimes\cdots\otimes\Delta^{d_r-1}$
par le morphisme d'Eilenberg-Zilber correspond \'egalement \`a l'image
du g\'en\'erateur $\Delta^{d_1-1}\otimes\cdots\otimes\Delta^{d_k-2}\otimes\cdots\otimes\Delta^{d_r-1}$
de $N_{d_1-1}(\Delta^{d_1-1})\otimes\cdots\otimes N_{d_k-2}(\Delta^{d_k-2})\otimes\cdots\otimes N_{d_r-1}(\Delta^{d_r-1})$
par l'application compos\'ee du morphisme d'Eilenberg-Zilber
$N_{d_1-1}(\Delta^{d_1-1})\otimes\cdots\otimes N_{d_k-2}(\Delta^{d_k-2})\otimes
\cdots\otimes N_{d_r-1}(\Delta^{d_r-1})
\,\rightarrow\,N_{d-1}(\Delta^{d_1-1}\times\cdots\times\Delta^{d_k-2}\times\cdots\times\Delta^{d_r-1})$
et du morphisme simplicial
$N_*(1\times\cdots\times d^x\times\cdots\times 1):
N_*(\Delta^{d_1-1}\times\cdots\times\Delta^{d_k-2}\times\cdots\times\Delta^{d_r-1})
\,\rightarrow\,N_*(\Delta^{d_1-1}\times\cdots\times\Delta^{d_k-1}\times\cdots\times\Delta^{d_r-1})$.
Une face $\Delta^{d_1-1}\otimes\cdots\otimes d_x(\Delta^{d_k-1})\otimes\cdots\otimes\Delta^{d_r-1}$
donne par suite un \'el\'ement de la forme $\TC(v)$ dans $N_*(\W(r))$:
on consid\`ere la surjection $v$
telle que $\tau_v = \tau_u\cdot 1\times\cdots\times d^x\times\cdots\times 1$
(cf. lemme \ref{SurjFaces}).
Les surjections $v$ ainsi obtenues sont les termes de la diff\'erentielle de $u$ dans le complexe des surjections $\X(r)$
(cf. [\bibref{BF}]).
On v\'erifie que la diff\'erence entre le signe de $v$ dans la diff\'erentielle de $u\in\X(r)$
et le signe de $\Delta^{d_1-1}\otimes\cdots\otimes d_x(\Delta^{d_k-1})\otimes\cdots\otimes\Delta^{d_r-1}$
dans la diff\'erentielle
de $\Delta^{d_1-1}\otimes\cdots\otimes\Delta^{d_r-1}\in N_*(\Delta^{d_1-1})\otimes\cdots\otimes N_*(\Delta^{d_r-1})$
correspond \`a la diff\'erence d'orientation entre les simplexes fondamentaux de $\tau_u$ et de $\tau_v$.
Ceci termine la preuve du lemme.

\tit{\sc Lemme}\label{AWPpty} --- {\it L'application $\TR: \E\,\rightarrow\,\X$ est caract\'eris\'ee
par les propri\'et\'es suivantes:

(1) c'est un morphisme de complexes;

(2) on suppose que $\sigma$ est le simplexe maximal d'un prisme $\tau_u$, alors:
$$\TR(\sigma) = \left\{\matrix{ u, & \ \hbox{si $\sigma$ est le simplexe fondamental de $\tau_u$},\hfill\cr
0, & \ \hbox{sinon}.\hfill\cr }\right.$$}

\smallskip
L'essentiel de la preuve de ce lemme est contenue dans l'article [\bibref{BF}].
L'ant\'ec\'edent d'une surjection $u$ que l'on construit dans cet article est en fait le simplexe fondamental de $\tau_u$.
On montre que le morphisme $\TR$ s'annule sur les autres simplexes maximaux de $\tau_u$
en reprenant nos arguments pour le calcul de l'image du simplexe fondamental.

Ces propri\'et\'es caract\'erisent l'application $\TR: \E\,\rightarrow\,\X$
parce que tout simplexe de $\W(r)$ est contenu dans un prisme $\tau_u$.
On en conclut \'egalement que le morphisme $\TR: \E\,\rightarrow\,\X$
est une application d'Alexander-Whitney
tordue par le choix du simplexe fondamental.

\references{\parindent=0cm\leftskip=1cm\rightskip=1cm

\biblabel{B}\refto{C. Berger},
{\it Op\'erades cellulaires et espaces de lacets it\'er\'es},
Ann. Inst. Fourier {\bf 46} (1996), 1125-1157.

\biblabel{BF} \refto{C. Berger, B. Fresse},
{\it Combinatorial operad actions on cochains},
pr\'e\-pu\-bli\-ca\-tion (2001).

\biblabel{GZ}\refto{P. Gabriel, M. Zisman},
Calculus of fractions and homotopy theory,
Ergebnisse der Mathematik und ihrer Grenzgebiete {\bf 35}, Springer-Verlag, 1967.

\biblabel{KS}\refto{M. Kontsevich, Y. Soibelman},
{\it Deformations of algebras over operads and De\-li\-gne's conjecture},
in ``Conf\'erence Mosh\'e Flato 1999, Vol. I'',
Math. Phys. Stud. {\bf 21}, Kluwer (2000), 255-307.

\biblabel{MCSPrism}\refto{J. McClure, J.H. Smith},
{\it A solution of Deligne's conjecture},
pr\'e\-pu\-bli\-ca\-tion (1999).

\biblabel{MCSSeq}\reftosame,
{\it Multivariable cochain operations and little $n$-cubes},
pr\'e\-pu\-bli\-ca\-tion (2001).

}

\bye